\documentclass[12pt,psfig, reqno]{amsart}
\usepackage{amsmath, amsthm, amsfonts, amssymb, latexsym, graphicx}
\usepackage[latin1]{inputenc}

\newtheorem{theorem}{Theorem}[section]
\newtheorem{lemma}[theorem]{Lemma}

\theoremstyle{definition}

\newtheorem{example}[theorem]{Example}

\numberwithin{equation}{section}

\newcommand{\mC}{\ensuremath{\mathbb{C}}}

\newcommand{\mN}{\ensuremath{\mathbb{N}}}

\newcommand{\mD}{\ensuremath{\mathbb{D}}}

 \allowdisplaybreaks

\author{}
\date{}

\begin{document}

\title{Normality from one family of meromorphic functions to another through sharing of values}
\author[K. S. Charak]{Kuldeep Singh Charak}
\address{
\begin{tabular}{lll}
&Kuldeep Singh Charak\\
&Department of Mathematics\\
&University of Jammu\\
&Jammu-180 006\\ 
&India\\
\end{tabular}}
\email{kscharak7@rediffmail.com}

\author[M. Kumar]{Manish Kumar}
\address{
\begin{tabular}{lll}
&Manish Kumar\\
&Department of Mathematics\\
&University of Jammu\\
&Jammu-180 006\\ 
&India\\
\end{tabular}}
\email{manishbarmaan@gmail.com}

\author[R. Kumar]{Rahul Kumar}
\address{
\begin{tabular}{lll}
&Rahul Kumar\\
&Department of Mathematics\\
&University of Jammu\\
&Jammu-180 006\\
&India
\end{tabular}}
\email{rktkp5@gmail.com}

\maketitle

\begin{abstract} Let $\mathcal F$ and $\mathcal G$ be two families of meromorphic functions on a domain $D$, and let $a,\ b$ and $c$ be three distinct points in the extended complex plane. Let $\mathcal G$ be a normal family in $D$ such that all limit functions of $\mathcal G$ are non-constant. If for each $f\in\mathcal F$, there exists $g\in\mathcal G$ such that  $f$ and $g$ share $a,\ b$ and $c$ partially, then $\mathcal F$ is normal in $D.$ This  gives a sharp improvement of a result due to X. J. Liu, S. H. Li and X. C. Pang. We also prove some interesting related sharp results.
\end{abstract}


\renewcommand{\thefootnote}{\fnsymbol{footnote}}

\footnotetext{The work of the first author is partially supported by Mathematical Research Impact Centric Support (MATRICS) grant, File No. MTR/2018/000446, by the Science and Engineering Research Board (SERB), Department of Science and Technology (DST), Government of India.}

\section{Introduction and Main Results.}

For the sake of brevity we shall use the following standard notations:
\begin{itemize}
	\item $\mathcal{H}(D):$ class of all holomorphic functions on a domain $D$ in $\mC;$
	\item $\mathcal{M}(D):$ class of all meromorphic functions on a domain $D$ in $\mC;$
	\item $E(f,a):$ the set of $a-$points of $f$ counted with multiplicity (CM);
	\item $\bar{E}(f,a):$ the set of $a-$points of $f$ counted by ignoring multiplicity (IM).
\end{itemize}

The present paper deals with passing of normality from one subfamily of $\mathcal{M}(D)$  to another when the two subfamilies satisfy a certain condition on sharing of values.

\smallskip

 Let's recall that a subfamily $\mathcal F$ of $\mathcal{M}(D)$ is said to be normal in $D$ if from any given sequence in $\mathcal F$ we can extract a subsequence that converges spherically locally uniformly in $D.$ The limit function is either in $\mathcal{M}(D)$ or identically equal to $\infty$. For complete introduction to normal families of meromorphic functions, the reader may refer to \cite{schiff} and \cite{zal}.

\smallskip

Two functions $f,g \in \mathcal{M}(D)$ are said to share a given value $a\in \mC_{\infty}$ IM (CM) if $\bar{E}(f,a)=\bar{E}(g,a)$ ($E(f,a)=E(g,a));$ further, $f$ and $g$ are said to share $a$ partially if  $\bar{E}(f,a)\subset \bar{E}(g,a).$ Maximum number of values shared by two non-constant and distinct meromorphic functions is $4$ (see \cite{H2}). For detailed account of sharing of values by meromorphic functions and their uniqueness, one may refer to \cite{yang}.

   A study of normality of two families of meromorphic functions on a common domain  was initiated in $2013$, almost simultaneously, by  Liu, Li and Pang \cite{Liu} and Yuan, Xiong and Lin \cite{yuan}. Liu, Li and Pang \cite{Liu} proved the following result:

\smallskip

\noindent
\textbf{Theorem A} {\it  Let $\mathcal F$ and $\mathcal G$ be two subfamilies of $\mathcal{M}(D)$ and let $a,\ b,\ c, \ d \in \mC$ be four distinct values. If $\mathcal G$ is normal in $D$, and for each $f\in\mathcal F$, there exists $g\in\mathcal G$ such that $f$ and $g$ share the values $a,\ b,\ c$ and $ d$, then $\mathcal F$ is normal in $D$.}

\smallskip

Though the number $4$ of shared values in Theorem A cannot be reduced to $3$,  we can compensate the sharing of some of the values by a condition on limit functions of $\mathcal{G}$ or on $a-$points of $\mathcal{F}$ . Precisely we have obtained:

\begin{theorem} \label{th6} Let $\mathcal F$ and $\mathcal G$ be two  subfamilies of $\mathcal{M}(D)$ and let $a_1,\ a_2, \ a_3\in\mC_{\infty}$ be three distinct values. Let $\mathcal G$ be  normal in $D$ such that all limit functions of $\mathcal G$ are non-constant. If for each $f\in\mathcal F$, there exists $g\in\mathcal G$ such that $f$ and $g$ share $a_j, \ j=1,2,3,$ partially in $D,$ then $\mathcal F$ is normal in $D$.
\end{theorem}

The condition that ``all the limit functions of $\mathcal G$ are non-constant" in Theorem \ref{th6} is essential as shown by the following example:

\begin{example}\label{ex1} Let $\mathcal F=\left\{f_n:n\in\mathbb N\right\}$, where $f_n(z)=\tan nz$, be the family of meromorphic functions on the open unit disk $\mathbb D$. Let $z_{n,1},\ldots, z_{n,m_n} $ be the zeros of $\tan nz$ in $\mathbb D$. Let $\mathcal G=\left\{g_n\right\}$, where 
$$g_n(z)=\frac{1}{n}\prod_{i=1}^{m_n}\left(\frac{z-z_{n,i}}{1-\bar{z}_{n,i}z}\right).$$
 Then $f_n$ omits $i,-i$ on $\mathbb D$ and  each $f_n$ and  $g_n$ share $0$ partially. It is easy to see that all subsequences of $\left\{g_n\right\}$ converge locally uniformly to 0. Since $|g_n(z)|< 1$, $\mathcal G$ is normal in $\mathbb D$. But $\mathcal F$ is not normal at 0. \end{example}

\begin{theorem}\label{thm1} Let $a_1,\ a_2, \ a_3, \ a_4 \in\mC_{\infty}$ be four distinct values, $\mathcal F$ and $\mathcal G$ be two subfamilies of $\mathcal{M}(D)$ such that the multiplicity of $a_1-$points of each $f\in\mathcal F$ is at least $2,$ and  let $\mathcal G$ be normal in $D$.  If for each $f\in\mathcal F$, there exists $g\in\mathcal G$ such that $f$ and $g$ share $a_j, \ j=2,3,4,$ partially in $D$, then $\mathcal F$ is normal in $D$.
\end{theorem}

Example \ref{ex1} also shows that in Theorem \ref{thm1}, the condition, ``the multiplicity of $a-points$ of $f\in\mathcal F$ is at least 2" is essential.

\begin{theorem}\label{thm3} Let $a_1,\ a_2, \ a_3, \ a_4 \in\mC_{\infty}$ be four distinct values, $\mathcal F$ and $\mathcal G$ be two subfamilies of $\mathcal{M}(D)$ such that the multiplicities of $a_1-$points and $a_2-$points of each $f\in\mathcal F$ are at least $2$ and $3$ respectively, and  let $\mathcal G$ be normal in $D$. If for each $f\in\mathcal F$, there exists $g\in\mathcal G$ such that $f$ and $g$ share $a_3$ and $a_4$ partially in $D$, then $\mathcal F$ is normal in $D$.
\end{theorem}

To pass on the normality of $\mathcal{F}\subseteq \mathcal{M}(D)$ to a family $\mathcal{R}$ of rational functions, we only require three values to be shared partially:

\begin{theorem}\label{thm2}
Let $\mathcal{R}$ be a family of rational functions of degree at most $m\in \mN$ and let  $\mathcal{G}\subseteq \mathcal{M}(D).$ Let $a_1,\ a_2, \ a_3\in \mC_{\infty}$ be three distinct  values. If $\mathcal{G}$ is normal in $D$ and for each $R\in\mathcal{R}$, there exists $g\in\mathcal{G}$ such that $R$ and $g$ share $a_j, \ j=1,2, 3,$ partially in $D$, then $\mathcal{F}$ is normal in $D.$
\end{theorem}

\begin{example}
Let $\mathcal{F}=\{f_n(z)=nz: n\in\mathbb{N}\} ~\mbox{and}~ \mathcal{G}=\{g_n(z)=z^n: n\in\mathbb{N}\}$. Then $\mathcal{F}$ is a family of rational map of degree $1.$ Also $\mathcal{G}$ is normal in $\mathbb{D}$ and for each $f_n$ there exists $g_n$ such that $f_n$ and $g_n$ share $0$ and $\infty$ partially.  But $\mathcal{F}$ is not normal in $\mathbb{D}.$ This shows that the number of values cannot be reduced to two. 
\end{example}

Liu, Li and Pang \cite{Liu} also proved the following two results on the  normality of two families: 

\noindent\textbf{Theorem B} {\it Let $\mathcal F, \ \mathcal G \subset \mathcal{H}(D)$, all of whose zeros have multiplicity at least $k+1, \ k\in \mN.$  Let $b$ be a non-zero complex number. Assume that $\mathcal G$ is normal, and for any subsequence $\left\{g_n\right\} \subseteq \mathcal G$ such that $g_n\to g$ locally uniformly on $D$, we have $g\not\equiv \infty$ and $g^{(k)}\not\equiv b$ on $D$. If for every $f\in\mathcal F$, there exists $g\in\mathcal G$ such that $f$ and $g$ share $0$, and $f^{(k)}$ and $g^{(k)}$ share $b$, then $\mathcal F$ is normal in $D$.}

\noindent\textbf{Theorem C} {\it Let $\mathcal F, \ \mathcal G \subset \mathcal{M}(D)$, all of whose zeros have multiplicity at least $k+1, \ k\in \mN,$ and let $b$ be a non-zero complex number. Suppose that $\mathcal G$ is normal in $D$ and for any sequence $\left\{g_n\right\}$ of $\mathcal G$ converging to a function $g$ spherically locally uniformly in $D$, we have $g^{(k)}\not\equiv b$ and $g\not\equiv\infty$. If for every $f\in\mathcal F$, there exists $g\in\mathcal G$ such that $f$ and $g$ share $0$ and $\infty$ IM, and  $f^{(k)}$ and $g^{(k)}$ share $b$ CM, then $\mathcal F$ is normal in $D$.}

\smallskip

We have obtained the following improvement of Theorem C:

\begin{theorem} \label{th2} Let $\mathcal F, \ \mathcal G \subset \mathcal{M}(D)$, all of whose zeros have multiplicity at least $k+1, \ k\in \mN,$ and let $a$ be a non-zero complex number.  Suppose that $\mathcal G$ is normal in $D$ such that all limit functions of $\mathcal G$ are not identically equal to infinity. If for each $f\in\mathcal F$, there exists $g\in\mathcal G$ such that $f$ and $g$ share $0$ and $\infty$ partially, and $f^{(k)}$ and $g^{(k)}$ share $a$ partially, then $\mathcal F$ is normal in $D$.
\end{theorem}

\noindent\textbf{Note:} When $\mathcal F$ and $ \mathcal G $ in Theorem \ref{th2} are restricted to be in $\mathcal{H}(D)$, we immediately obtain an improvement of Theorem B.

\smallskip

 One can show that condition, ``partial sharing of $0$ and $\infty$" in Theorem \ref{th2}, is essential, for example one may refer to Examples 1.7 and 1.8  in \cite{Liu}.

\begin{example} Consider 
$$\mathcal F :=\left\{f_n(z)=e^{nz}, \ n\in\mathbb N, \ z\in \mathbb D \right\}$$
and  for $k\in \mN$, consider
$$\mathcal G:=\left\{g_n(z)=(z-2)^{(k+1)}, \ n\in\mathbb N, z\in \mathbb D \right\}.$$
Then each $f\in\mathcal F$ omits 0 and $\infty$, $\mathcal G$ is normal in $\mathbb D$ and all limit functions of $\mathcal G$ are not identically equal to infinity. But $\mathcal F$ is not normal in $\mathbb D$. This shows that the condition, ``$f^{(k)}$ and $g^{(k)}$ share $a$, partially" in Theorem \ref{th2} is essential.
\end{example}
 
\begin{example} Consider the family 
$$\mathcal F:=\left\{f_n(z)=z-\frac{e^{nz}}{n}, \ n\in \mN, \ z\in \mD\right\}$$
and let $\left\{z_{n,i}\right\}_{i=1}^{k_n}$ be the zeros of $f_n$ in $\mathbb D$. Further consider 
$$\mathcal G:=\left\{g_n(z)= \prod^{k_n}_{i=1}\left(\frac{z-z_{n,i}}{1-\overline{z_{n,i}}z}\right)^{2}, n\in \mN, z\in \mD \right\}.$$
Then $\mathcal F, \mathcal G \subseteq \mathcal{H}(\mD)$ such that all the zeros of $f\in\mathcal F$ are simple and all the zeros of $g\in\mathcal G$ have multiplicity $2$. It is easy to see that for each $f_n$, there exists $g_n\in\mathcal G$ such that $f_n$ and $ g_n$ share $0$ partially and $f^{'}_n$ and $g^{'}_n$ share $1$ partially. Since $|g_n(z)|<1$, $\mathcal G$ is normal in $\mathbb D$ and all  limit functions of $\mathcal G$ are not identically equal to infinity. But $\mathcal F$ is not normal in $\mathbb D$. This shows that the condition, ``all the zeros of $f\in\mathcal F$ have multiplicity at least $k+1$"  in Theorem \ref{th2}, is essential.
\end{example}

\begin{example} Let $\mathcal F:=\left\{f_n(z)=nz^{2k}\right\}$ and  $\mathcal G:=\left\{g_n(z)=z^{k}/k!\right\}.$
Then $\mathcal F, \mathcal G \subseteq \mathcal{H}(\mD)$  with all zeros of $f\in\mathcal F$ of multiplicity at least $k+1$ and all  zeros of $g$ have multiplicity exactly $k$. It is easy to see that for each $f\in\mathcal F$, there exists $g\in\mathcal G$ such that $f$ and $g$ share $0$ partially,  and $f^{(k)}$ and $g^{(k)}$ share $1$ partially. Also $\mathcal G$ is normal in $\mathbb D$ and all limit functions of $\mathcal G$ are not identically equal to infinity. But $\mathcal F$ is not normal in $\mathbb D$. This shows that in Theorem \ref{th2} the condition, ``all the zeros of $g\in\mathcal G$ have multiplicity at least $k+1$", is essential.
\end{example}

Also,  the condition, ``all the limit functions of $\mathcal G$ are not identically equally to infinity" in Theorem \ref{th2} is essential:

\begin{example} Let $\mathbb D$ be the open unit disk and let $k$ be a positive integer. Let 
$$\mathcal F:=\left\{\frac{e^{nz}}{e^{nz}-1}:n\in\mathbb N\right\}$$ and 
$$\mathcal G:=\left\{\frac{e^{nz}}{e^{nz}-1}+n:n\in\mathbb N\right\}$$
 be two families of meromorphic functions on $\mathbb D$. Clearly, it satisfies all the conditions of Theorem \ref{th2} except ``all the limit functions of $\mathcal G$ are not identically equal  to infinity". It is easy to see that $\mathcal F$ is not normal at $0$. 
\end{example}

Finally, we have the following related results:

\begin{theorem}\label{th3} Let $\mathcal F$ and $\mathcal G$ be two subfamilies of $\mathcal{M}(D)$ all of whose zeros have multiplicity at least $k+1$. Let $a_1$ and $a_2$ be two complex numbers, and let $\mathcal G$ be normal in $D$ such that all limit functions of $\mathcal G$ are not identically equal to infinity. If, for each $f\in\mathcal F$, there exists $g\in\mathcal G$ such that $f$ and $g$ share $\infty$ partially,  and $f^{(k)}$ and $g^{(k)}$ share $a_j, \ j=1,2,$ partially, then $\mathcal F$ is normal in $D$.
\end{theorem}

\begin{theorem} \label{th4} Let $\mathcal F \subseteq \mathcal{M}(D)$ be such that all its zeros have multiplicity at least $k+1, \ k\in \mN$, and let  $a_1, a_2$ and $a_3$ be three distinct complex numbers. Let $\mathcal G \subset \mathcal{M}(D)$ be a normal family in $D$ such that all limit functions of $\mathcal G$ are not identically equal to infinity. If for each $f\in\mathcal F$, there exists $g\in\mathcal G$ such that 
$f^{(k)}$ and $g^{(k)}$ share $a_j, \ j=1,2,3,$  partially, then $\mathcal F$ is normal in $D$.
\end{theorem}

Example \ref{ex1} also shows that the condition, ``all the limit functions of $\mathcal G$ are not identically equal to infinity" in Theorem \ref{th3} and Theorem \ref{th4}, is essential.

\begin{example}  Let 
$$\mathcal F:=\left\{f_n(z)=nz^{k}:n\in\mathbb N, \ z\in \mD \right\}$$ and let 
$$\mathcal G:=\left\{g_n(z)=z+1/n: n\in\mathbb N, \ z\in \mD\right\}.$$
Then  $f_n^{(k)}$ and $g_n^{(k)}$ share any three given values partially, outside the set $\{k!n:n\in\mathbb N\}$. Also, $\mathcal G$ is normal and all the limit functions of $\mathcal G$ are not identically equal to infinity. But $\mathcal F$ is not normal at $0$. This shows that in Theorem \ref{th4} the condition, ``all the zeros of functions in $\mathcal F$ have multiplicity at least $k+1$" is essential.
\end{example}
The condition, ``$f$ and $g$ share $\infty$ partially" in Theorem \ref{th3} is essential. Also, the number of shared values in Theorem \ref{th4} cannot be reduced to two, for example one may refer to \cite{Liu}.

\section{Proofs of Main Theorems}

Besides Zalcman's Lemma \cite{zal} we shall use the following results-stated as lemmas-in the proofs of our main theorems:

\begin{lemma} \label{l4} \cite{W3} Let $f$ be a transcendental meromorphic function of finite order on $\mathbb C$ all of whose zeros have multiplicity at least $k+1$, where $k$ is a positive integer. Then $f^{(k)}$ assumes every non-zero complex number infinitely many times on $\mathbb C$.
\end{lemma}

\begin{lemma}\label{thmN}\cite{Bergweiler}
Let $a_1, \ldots, a_q\in \mC_\infty \mbox{ and } m_1, \ldots, m_q \in \mN$, where $q$ is a positive integer. Suppose $f\in \mathcal{M}(\mC)$ is non-constant such that all $a_j$-points of $f$ have multiplicity  at least $m_j, (j=1, \ldots, q).$ Then 
$$\sum\limits_{j=1}^{q}\left(1-\frac{1}{m_j}\right) \leq 2.$$
\end{lemma}

\noindent If $f$ does not assume the value $a_j$, then we take $m_j=\infty$.

\begin{lemma}\cite{pang} \label{l1} Let $\mathcal{F}$ be a family of meromorphic functions in the open unit disk $\mathbb{D}$ with the property that for each $f\in \mathcal F$, all zeros of $f$ are of multiplicity at least $k$. Suppose that there exists a number $A\geq 1$ such that $|f^{k}(z)|\leq A$ whenever  $f\in \mathcal{F}$ and $f(z)=0$. If $\mathcal{F}$ is not normal in $\mD$, then there exist, for each $0\leq \alpha \leq k$, a number $0<r<1$, points $z_n$ with $|z_n|<1$, functions $f_n\in \mathcal{F}$ and positive numbers $\rho_n\rightarrow 0$ such that 
$$g_n(\zeta)=\rho_{n}^{-\alpha}f_n(z_n+\rho_n\zeta)\rightarrow g(\zeta)$$
  spherically locally uniformly on $\mathbb{C}$, where $g$ is a non-constant meromorphic function on $\mathbb C$.
\end{lemma}

\begin{lemma}\label{l5}\cite{wang} Let $f$ be a non-constant meromorphic function of finite order on $\mathbb C$, all of whose zeros have multiplicity at least $k+1$. If $f^{(k)}(z)\neq a$ on $\mathbb C$, where $a\in\mathbb C\setminus \left\{0\right\}$, then 
$$f(z)=\frac{a\ (z-b)^{k+1}}{k!\ (z-c)},$$
where $b$ and $c$ are two distinct complex numbers.
\end{lemma}

\medskip

 \textbf{Proof of Theorem \ref{th6}} Suppose that $\mathcal F$ is not normal at $z_0\in D$. Then, by Zalcman Lemma, there exist $z_n\to z_0$, $f_n\in\mathcal F$ and $\rho_n\to 0^{+}$ such that $F_n(\zeta)=f_n(z_n+\rho_n\zeta)$ converges spherically locally uniformly to a non-constant meromorphic function $F$ in $\mathbb C$. Now by hypothesis, there exists a sequence $\{g_n\}\subset \mathcal G$  such that for each $n,$ $f_n$ and $g_n$ share $a_j, \ j=1,2,3$, partially in $D$. We may assume that $\{g_n\}$ converges spherically locally uniformly to a non-constant meromorphic function $g$ in $D$.

\smallskip

 {\it Claim.} $F$ assumes at most one of the values $a_j, \ j=1,2,3.$ 

\noindent Supposing on the contrary, we assume that there exist $\zeta_j\in \mathbb C$ such that $F(\zeta_j)=a_j$ for $j=1,2.$ Then Hurwitz's theorem ensures the existence of sequences $\{\zeta_{j,n}\}:\zeta_{j,n}\to\zeta_j$ such that for sufficiently large $n,$
$F_n(\zeta_{j,n})=a_j, j=1,2.$ Thus, by hypothesis, we find that $g_n(z_n+\rho_n\zeta_{j,n})=a_j, j=1,2.$
That is, $g(z_0)=a_j, \ j=1,2$ which is not allowed and hence the claim.

\smallskip

We assume that $F$ omits $a_1$ and $a_2$. Then $F$ is  transcendental and $F(\zeta_0)=a_3$, for some $\zeta_0\in \mathbb C$. Using the preceding argument, we arrive at $g(z_0)=a_3$. Let $m\geq 1$ be the order of zero of $g-a_3$ at $z_0$. Then again by Hurwitz's theorem there are exactly $m$ zeros of $g_n-a_3$ in a neighborhood $N(z_0)$ of $z_0,$ for sufficiently large $n.$  Since $F$ assumes $a_3$ infinitely often, let $\zeta_j, j=1,\ldots, m+1$, be the distinct complex numbers such that $F(\zeta_j)=a_3$. Once again by Hurwitz's theorem, for each $j$ there exist  $\{\zeta_{j,n}\}$ such that $\zeta_{j,n}\to\zeta_j$ and  $F_n(\zeta_{j,n})=a_3$ for sufficiently large $n.$ That is, $f_n-a_3$  and hence $g_n-a_3$ has $m+1$ distinct zeros in $N(z_0)$, namely, at $z_n+\rho_n\zeta_{j,n}, j=1,\ldots, m+1$. This is a contradiction. $\Box$

\medskip

Proofs of Theorem \ref{thm1}, Theorem \ref{thm3} and Theorem \ref{thm2}: Suppose that $\mathcal F$ is not normal in the domain $D$. By applying the arguments of the proof of the Theorem \ref{th6} we find that limit function $F$ omits two values in Theorem \ref{thm1}, one value in Theorem \ref{thm3} and two values in Theorem \ref{thm2}. Then by Argument Principle and Lemma \ref{thmN} we arrive at a contradiction in Theorem \ref{thm1} and Theorem \ref{thm3}, whereas in Theorem \ref{thm2} we arrive at a contradiction by using Hurwitz's theorem only.
 
\medskip

\textbf{Proof of the Theorem \ref{th2}}  Since  normality is a local property, we assume $D$ to be the open unit disk $\mD.$ Suppose that $\mathcal F$ is not normal at $z_0\in \mD$. Then, by Lemma \ref{l1}, there exist $z_n\to z_0$, $f_n\in\mathcal F$ and $\rho_n\to 0^{+}$ such that $F_n(\zeta)=\rho_n^{-k}f_n(z_n+\rho_n\zeta)$ converges spherically locally uniformly to $F$ in $\mC$, where $F$ is a non-constant meromorphic function of finite order on $\mathbb C$ and all the zeros of $F$ have multiplicity at least $k+1$.\\
Without loss of generality, we may assume that the corresponding sequence $\left\{g_n\right\}\subset\mathcal G$ converges spherically locally uniformly to $g$ in $\mD$, where $g\not\equiv\infty$ and  all the zeros of $g$ have multiplicity at least $k+1$.

\smallskip

{\it Case-$1:$} When $F$ assumes 0.\\
Suppose that $F(\zeta_0)=0$ for some $\zeta_0\in\mathbb C$. Then, by Hurwitz's theorem, there exists $\zeta_n\to\zeta_0$ such that for sufficiently large $n,$ $F_n(\zeta_n)=0.$ That is, $f_n(z_n+\rho_n\zeta_n)=0$, and hence $g_n(z_n+\rho_n\zeta_n)=0,$ for sufficiently large $n.$ That is, $g(z_0)=0.$ Since  zeros of $g$ have multiplicity at least $k+1$, $g^{(k)}(z_0)=0.$\\
{\it Claim: } $F^{(k)}(\zeta)\neq a$ on $\mathbb C$.\\
Suppose that there exists $\zeta_1\in\mathbb C$ such that $F^{(k)}(\zeta_1)=a.$ If $F^{(k)}\equiv a$, then $F$ is a polynomial of degree $k$, a contradiction to the fact that all the zeros of $F$ are of multiplicity at least $k+1$. Therefore, $F^{(k)}\not\equiv a$, and so by Hurwitz's theorem, there exists $\zeta_{1,n}:\zeta_{1,n}\to\zeta_1$  such that for sufficiently large $n$, $F^{(k)}_n(\zeta_{1,n})=a.$  This  by our hypothesis implies that  $g^{(k)}_n(z_n+\rho_n\zeta_{1,n})=a$, and so $g^{(k)}(z_0)=a \ (\neq 0)$, a contradiction. This proves the claim.\\
By Lemma \ref{l5}, $F$ is a rational function and  therefore, there exists $\zeta_2\in\mathbb C$ such that $F(\zeta_2)=\infty$. Again by  Hurwitz's theorem, there exists $\{\zeta_{2,n}\}:\zeta_{2,n}\to\zeta_2$  such that for sufficiently large $n$, $F_n(\zeta_{2,n})=\infty.$ This by our hypothesis implies that for sufficiently large $n,$ $g_n(z_n+\rho_n\zeta_{2,n})=\infty,$ and hence $g(z_0)=\infty$, a contradiction.

\smallskip

{\it Case-$2:$} When $F$ omits zero.\\
 In this case we first prove that $F(\zeta)\neq \infty$ on $\mathbb C.$ For,  suppose that $F(\zeta_3)=\infty$, for some $\zeta_3\in\mathbb C.$ Then there exists $\{\zeta_{3,n}\} : \zeta_{3,n}\to\zeta_3$ such that $F_n(\zeta_{3,n})=\infty  $, for sufficiently large $n$. This by our hypothesis implies that $g_n(z_n+\rho_n\zeta_{3,n})=\infty$, and so $g(z_0)=\infty$. Since $F\neq 0$ on $\mathbb C$, by Hayman's Alternative\cite{H1}, $F^{(k)}(\zeta_4)=a$, for some $\zeta_4\in\mathbb C$. By Hurwitz's theorem, there exists a sequence $\{\zeta_{4,n}\}$ converging to $\zeta_4$ such that for sufficiently large $n,$  $F_n^{(k)}(\zeta_{4,n})=a,$  and hence by hypothesis we have $g_n^{(k)}(z_n+\rho_n\zeta_{4,n})=a$, and so $g^{(k)}(z_0)=a$, a contradiction. Therefore, $F$ is transcendental and again by Hayman's Alternative,  $F^{(k)}$ assumes $a$. As in Case-$1$, we get $g^{(k)}(z_0)=a$. If $g^{(k)}\equiv a$, then $g$ is a polynomial of degree $k$, a contradiction to the fact that zeros of $g$ are of multiplicity at least $k+1$. Assume that $m\ (\geq 1)$ be the order of zero of $g^{(k)}-a$ at $z_0$. By Hurwitz's theorem, for sufficiently large $n$, $g_n^{(k)}-a$ has exactly $m$ zeros in a neighborhood $N(z_0)$ of $z_0$.  Now by Lemma \ref{l4}, $F^{(k)}$ assumes $a$ infinitely many times and let $\zeta^{'}_1,\ldots, \zeta^{'}_{m+1}$ be the distinct zeros of $F^{(k)}-a$. By Hurwitz's theorem, there exists $\{\zeta^{'}_{j,n}\}$ converging to $\zeta^{'}_j$ such that for sufficiently large $n,$ $F_n^{(k)}(\zeta^{'}_{j,n})=b,\ j=1,\ldots, m+1.$
This implies that $f^{(k)}_n-a$ has $m+1$ distinct zeros in $N(z_0)$, and hence $g^{(k)}_n-a$ has $m+1$ distinct zeros in $N(z_0),$ a contradiction. $\Box$

\bigskip

\textbf{Proof of the Theorem \ref{th3}} Normality being a local property allows to replace the domain $D$ by the open unit disk $\mD.$ By supposing on the contrary that $\mathcal{F}$ is not normal at $z_0 \in \mD,$  Lemma \ref{l1} ensures the existence of sequences $\{z_n\}$ in $\mD,$   $\{f_n\}$ in $\mathcal F$ and $\{\rho_n\}$ in $(0,1)$ such that $z_n \to z_0, \ \rho_n \to 0^{+}$ and  
$$F_n(\zeta):=\frac{f_n(z_n+\rho_n\zeta)}{\rho_{n}^{k}}\to F(\zeta),$$
spherically locally uniformly on $\mathbb C$, where $F$ is a non-constant meromorphic function on $\mathbb C$ all of whose zeros are of multiplicity at least $k+1$.

We may assume that the corresponding sequence $\{g_n\}$ converges  spherically locally uniformly to $g\not\equiv\infty$ on $\mD$. For the sake of convenience we denote the cardinality of a set $A$ by $\#(A).$

\smallskip

{\it Claim:} $\#\left(F^{(k)}(\mathbb C)\cap \left\{a_1, a_2\right\}\right)\leq 1$.\\
Suppose that there exist $\zeta_j\in\mathbb C$ such that $F^{(k)}(\zeta_j)=a_j, \ j=1,2.$  If $F^{(k)}$ is a constant function, then $F$ is a polynomial of degree at most $k$ which is not the case as all the zeros of $F$ are of multiplicity at least $k+1$. Now  Hurwitz's theorem implies the existence of sequences $\{\zeta_{j,n}\}$ such that $F_n^{(k)}(\zeta_{j,n})=a_j, \ j=1,2,$ which in turn implies that $g^{(k)}_n(z_n+\rho_n\zeta_{j,n})=a_j$ and hence  $g^{(k)}(z_0)=a_j, \ j=1,2;$ which is absurd. This proves the claim.

\smallskip

 We may assume that $F^{(k)}$ omits $a_1$. Then  $F^{(k)}$ assumes $a_2,$ otherwise  $F^{(k)}$ would reduce to a constant enforcing $F$ to be a polynomial of degree at most $k$, a contradiction since all zeros of $F$ are of multiplicity at least $k+1$.  Thus there exists $\zeta_0\in\mathbb C$ such that $F^{(k)}(\zeta_0)=a_2$ which implies that $g^{(k)}(z_0)=a_2.$ Now we have two cases: First, when  $F$ assumes $\infty.$ In this case, since $f$ and $g$ share $\infty$ partially, $g(z_0)=\infty$ which is not true as  $g^{(k)}(z_0)=a_2$ and $g\not\equiv\infty.$ The second case when  $F$ omits $\infty.$ Then  $F^{(k)}$ assume $a_2$ infinitely many times, and one can arrive at a contradiction by the arguments used in the proof of Theorem \ref{th2}. $\Box$

\bigskip

\textbf{Proof of the Theorem \ref{th4}}  As in the last proof, we obtain a non-constant meromorphic function $F$ on $\mC$ all of whose zeros are of multiplicity at least $k+1$. Precisely $F$ is the spherical local uniform limit of the scaled sequence $F_n(\zeta):=\rho_{n}^{-k}f_n(z_n+\rho_n\zeta)$ in $\mD.$ Further, one can establish that 
 $$\#\left(F^{(k)}(\mathbb C)\cap \left\{a_j:j=1,2,3\right\}\right)\leq 1.$$

Thus $F^{(k)}$ omits two values in $\mathbb C$ and hence  $F^{(k)}$ reduces to a constant. This implies that $F$ is a polynomial of degree at most $k$, a contradiction as all zeros of $F$ are of multiplicity at least $k+1$. $\Box$

\end{document}